\documentclass{article}

\usepackage{amsfonts}
\begin{document}

\newcommand{\rf}[1]{(\ref{#1})}
\def\proof{\noindent {\bf Proof.}\ }

\def\weak{\mathring{\to}} 
\def\rr{\mathbb R} 
\def\rrr{{\bar \mathbb R}} 
\def\supp{\rm supp}
\def\te{\mathbb T} 

\def\endproof{\hfill \hspace{-6pt}\rule[-14pt]{6pt}{12pt}
\vskip2pt plus3pt minus 3pt}
\newtheorem{theorem}{Theorem}
\newtheorem{conjecture}[theorem]{Conjecture}
\newtheorem{lemma}[theorem]{Lemma}
\newtheorem{definition}[theorem]{Definition}
\newtheorem{definitions}[theorem]{Definitions}
\newtheorem{notation}[theorem]{Notation}
\newtheorem{notations}[theorem]{Notations}
\newtheorem{proposition}[theorem]{Proposition}
\newtheorem{remark}[theorem]{Remark}



\newcommand{\abs}[1]{\lvert#1\rvert}
\newcommand{\norm}[1]{\lVert#1\rVert}

\newcommand{\R}{\mathbb R}
\newcommand{\I}{\mathbb I}
\newcommand{\Rb}{\overline{\R}}
\newcommand{\C}{\mathbb C}
\newcommand{\Cb}{\overline{\C}}
\newcommand{\Z}{\mathbb Z}
\newcommand{\N}{\mathbb N}
\newcommand{\D}{\mathbb D}
\newcommand{\J}{{\mathcal J}}
\newcommand{\K}{{\mathcal K}}
\newcommand{\F}{{\mathcal F}}
\newcommand{\Rs}{{\mathcal R}}
\newcommand{\ie}{i.~e.,\ }
\renewcommand{\Re}{\operatorname{Re}}
\renewcommand{\Im}{\operatorname{Im}}

\setcounter{equation}{18}

\title{A Weierstrass-type theorem for homogeneous polynomials} \author{David Benko, Andr\'as Kro\'o$^\dagger$}
\maketitle

\begin{abstract}
 By the celebrated Weierstrass Theorem the set of algebraic polynomials
is dense in the space of continuous functions on a compact set in $\rr^d$. In this paper we study the following question: does the density hold if we approximate only by homogeneous polynomials? Since the set of homogeneous polynomials is nonlinear this leads to a nontrivial problem. It is easy to see that: 1) density may hold only on star-like {\bf 0}-symmetric
surfaces; 2) at least 2 homogeneous polynomials are needed for approximation. The most interesting special case of a star-like surface is a convex surface. It has been conjectured by the second author that functions continuous on {\bf 0}-symmetric convex surfaces in $\rr^d$ can be approximated by a pair of homogeneous polynomials. This conjecture is not resolved yet but we make substantial progress towards its positive settlement. In particular, it is shown in the present paper that the above conjecture holds for \ 
1) $d=2$, 2) convex surfaces in $\rr^d$ with $C^{1+\epsilon}$ boundary.
\end{abstract}

\footnotetext[1]
{$\dagger$ Supported by the OTKA grant \# T049196. 
Written partially during this author's visit at the
Center for Constructive Approximation, Vanderbilt University,
Nashville, Tennessee }

\footnotetext[2]
{AMS classification: Primary 41A10, 31A05, Secondary 52A10, 52A20.}

\footnotetext[3]
{Key words and phrases: Weierstrass, uniform approximation, homogeneous polynomials, convex body}

\section{Introduction}  The celebrated theorem of Weierstrass on
the density of real algebraic polynomials in the space of real
continuous functions on an interval $[a,b]$ is one of the main
results in analysis. Its generalization for real multivariate
polynomials was given by Picard, subsequently the
Stone-Weierstrass theorem led to the extension of these results
for subalgebras in $C(K)$.

In this paper we shall consider the question of density of
\emph{homogeneous polynomials}. Homogeneous polynomials are a
standard tool appearing in many areas of analysis, so the question
of their density in the space of continuous functions is a natural
problem. Clearly, the set of homogeneous polynomials is
substantially smaller relative to all algebraic polynomials. More
importantly, this set is nonlinear, so its density can not be
handled via the Stone-Weierstrass theorem. Furthermore, due to the
special structure of homogeneous polynomials some restrictions
should be made on the sets were we want to approximate (they have
to be star-like), and at least 2 polynomials are always needed for
approximation (an even and an odd one).

On the 5-th International Conference on Functional Analysis and
Approximation Theory (Maratea, Italy, 2004) the second author
proposed the following conjecture.

\begin{conjecture}\label{k0}  
Let $K\subset \R^{d}$ be a convex body which is centrally 
symmetric to the origin.   
Then for any function $f$ continuous on the
boundary $ Bd(K) $  
of $K$ and any $\epsilon >0$ there exist two
homogeneous polynomials $h$ and $g$ such that $|f-h-g|\leq
\epsilon$ on $ Bd(K) $.
\end{conjecture}

From now on we agree on the terminology that by ``centrally symmetric" we mean ``centrally symmetric to the origin".  

Subsequently in  \cite{ksz}  the authors verified the above Conjecture for
\emph{crosspolytopes} in $\R^{d}$ and arbitrary convex polygons in
$\R^{2}$.

In this paper we shall verify the Conjecture for those convex
bodies in $\R^{d}$ whose boundary $ Bd(K) $ is $C^{1+\epsilon}$ for
some $0<\epsilon\leq 1$ (Theorem \ref{k1}). Moreover, the Conjecture will
be verified in its full generality for $d=2$ (Theorem \ref{k2}).

It should be noted that parallel to our investigations P. Varj\'u
 \cite{varju}  also proved the Conjecture for $d=2$. In addition, he gives in
 \cite{varju}  an affirmative answer to the Conjecture for arbitrary centrally
symmetric polytopes in $\R^{d}$, and for those convex bodies in
$\R^{d}$ whose boundary is $C^{2}$ and has positive curvature. We
also would like to point out that our method of verifying the
Conjecture for $d=2$ is based on the potential theory and is
different from the approach taken in  \cite{varju}  (which is also based on
the potential theory). Likewise our method of treating
$C^{1+\epsilon}$ convex bodies is different from the approach used
in  \cite{varju}  for $C^{2}$ convex bodies with positive curvature.

\section{Main Results} 
Let $K$ be a centrally symmetric convex body
in $\R^{d}$. We may assume that $2\leq d$ and $\dim(K)=d$.  
The boundary of $K$ is $ Bd(K) $ which is given by the
representation $$ Bd(K) :=\{{\bf u} r({\bf u}): {\bf u}\in S^{d-1}\}$$ where $r$ is a
positive even real-valued function on $S^{d-1}$.
Here $S^{d-1}$
stands for the unit sphere in $\R^{d}$. We shall say that $K$ is
$C^{1+\epsilon}$, written $K \in C^{1+\epsilon}$, if the first   partial
derivatives of $r$ satisfy a Lip$\epsilon$ property on the unit
sphere, $\epsilon>0$. Furthermore denote by 
$$H^{d}_{n}:=
\Big\{\displaystyle\sum_{k_{1}+...+k_{d}=n}c_{\textbf{k}}\textbf{x}^{\textbf{k}}:
c_{\textbf{k}}\in \R\Big\}$$ 
the space of real homogeneous polynomials
of degree $n$ in $\R^{d}$. Our first main result is the following.

\begin{theorem}\label{k1}  
Let $K\in C^{1+\epsilon}$ be a centrally symmetric convex body in $\R^{d}$, where $0<\epsilon \leq 1$. Then
for every $f\in C( Bd(K) )$ there exist $h_{n}\in
H^{d}_{n}+H^{d}_{n-1}, n\in \N$ such that $ h_{n} \rightarrow f$
uniformly on $ Bd(K) $ as $n\rightarrow \infty.$ 
\end{theorem}

Thus Theorem \ref{k1} gives an affirmative answer to the Conjecture under
the additional condition of $C^{1+\epsilon}$ smoothness of the
convex surface. For $d=2$ we can verify the Conjecture in its full
generality. Thus we shall prove the following.

\begin{theorem}\label{k2}   
Let K be a centrally symmetric convex body
in $\R^{2}$. Then for every $f\in C( Bd(K) )$ there exist $h_{n}\in
H_{n}^{2}+H_{n-1}^{2},n\in \N$ such that $h_{n}\rightarrow f$
uniformly on $ Bd(K) $ as $n\rightarrow \infty.$ 
\end{theorem}

We shall see that Theorem \ref{k2} follows from

\begin{theorem}\label{k10}   
Let $1/W(x)$ be a positive convex function on $\rr$ such that 

\noindent
$|x|/W(-1/x)$ is also positive and convex. Let $g(x)$ be a continuous function which has the same limits at $-\infty$ and at $+\infty$. Then we can approximate $g(x)$ uniformly on $\rr$ by weighted polynomilas $W(x)^n p_n(x)$, $n=0,2,4,...$, $\deg p_n \le n$.
\end{theorem}

\section{Proof of Theorem \ref{k1}}
The proof of Theorem \ref{k1} will be based on several lemmas. The main
auxiliary result is the next lemma which provides an estimate for
the approximation of unity by even homogeneous polynomials. In
what follows $||...||_{D}$ stands for the uniform norm on $D$.

Our main lemma to prove Theorem \ref{k1} is the following.

\begin{lemma}\label{k3}  
Let $\tau \in (0,1)$. 
Under conditions of Theorem \ref{k1} there
exist $h_{2n}\in H^{d}_{2n}, n\in \N$, such that
$$||1-h_{2n}||_{ Bd(K) } = o(n^{-\tau\epsilon}).$$
\end{lemma}

The next lemma provides a partition of unity which we shall need
below. In what follows a cube in $\R^{d}$ is called regular if all
its edges are parallel to the coordinate axises. We denote the set 
$\{0,1,2,...\}^d$ by $\Z^{d}_{+}$.  

\begin{lemma}\label{k4}  
Given $0<h\leq 1$ there exist non-negative   
even functions $g_{\textbf{k}}\in C^{\infty}(\R^{d})$ such that their
support consists of $2^{d}$  
regular cubes with edge $h$, at most $2^{d}$
of supports of $g_{\textbf{k}}$'s have nonempty intersection, and
$$\displaystyle\sum_{\textbf{k}\in\Z^{d}_{+}}g_{\bf k}(\textbf{x})=1, \quad
\textbf{x}\in\R^{d}, \eqno(1)$$
$$|\partial^{m}g_{\textbf{k}}(\textbf{x})/\partial x_{j}^{m}|\leq c/h^{m},\quad {\bf x} \in \R^{d},  m\in \Z_{+}^{1},  1\leq j\leq d, \eqno(2)$$ where
$c>0$ depends only on $m \in \Z^{1}_{+}$ and $d$.
\end{lemma}

For the centrally symmetric convex body $K$ let
$$|\textbf{x}|_{K}:= \inf\{a>0: \textbf{x}/a \in K \}$$ be its Minkowski functional and set 
$$\delta_{K}:= \sup\{|\textbf{x}|/|\textbf{x}|_{K}: \textbf{x}\in
\R^{d}\} = \max\{|{\bf x} | : \ {\bf x}  \in Bd(K) \}.  $$ 
Moreover for $a\in  Bd(K) $ denote by $L_{a}$ a supporting
hyperplane at $a$.

 \begin{lemma}\label{k5}  
Let $\textbf{a}\in  Bd(K) , h_{n}\in H^{d}_{2n}$ be such that
 for any $\textbf{x}\in L_{\textbf{a}}, |\textbf{x}-\textbf{a}|\leq 4\delta_{K}$
 we have $|h_{n}(\textbf{x})|\leq 1$. Then whenever $\textbf{x}\in L_{\textbf{a}}$ satisfies
 $|\textbf{x}-\textbf{a}|>4\delta_{K}$ and $\textbf{x}/t \in K$ we have $$|h_{n}(\textbf{x}/t)|\leq (2/3)^{2n}. \eqno(3)$$
\end{lemma}

\begin{lemma}\label{k6}  
Consider the functions $g_{\textbf{k}}$
 from Lemma \ref{k4}. Then for at most $8^{d}/2h^{d}$ of them their support has
 nonempty intersection with $S^{d-1}$.
\end{lemma}

We shall verify first the technical Lemmas \ref{k4}-\ref{k6}, then the proof of Lemma \ref{k3}  will be given. 
Finally it will be shown that
Theorem \ref{k1} follows easily from Lemma \ref{k3}.

\textbf{Proof of Lemma \ref{k4}.} The main step of the proof consists of
verifying the lemma for $d=1$. Let $g\in C^{\infty}(\R)$ be an odd
function on $\R$ such that $g=1$ for $x<-1/2$ and monotone
decreasing from 1 to 0 on $(-1/2,0)$. Further, let $g^{*}(x)$ be
an even function on $\R$ such that $g^{*}(x)$ equals 1 on [0,1],
$g(x-3/2)/4+3/4$ on [1,2], and $g(x-5/2)/4+1/4$ on [2,3]. Then it
is easy to see that $g^{*}\in C^{\infty}(\R)$, it equals 1 on
$[-1,1]$, 0 for $|x|>3$ and is monotone decreasing on [1,3].
Moreover $$g^{*}(x) + g^{*}(x-4) = 1,\quad x\in [-1,5]. \eqno(4)$$
Set now $$g_{k}(x):= g^{*}(x-4k)+g^{*}(x+4k),\quad  k\in
\Z^{1}_{+}.$$ Then $g_{k}$'s are even functions which by (4)
satisfy relation
$$\displaystyle\sum_{k=0}^{\infty}g_{k}(x) = 1, \quad x\in \R.$$ In
addition, the support of $g_{k}$ equals $\pm [-3+4k,3+4k]$ and at
most 2 of $g_{k}$'s can be nonzero at any given $x\in \R$.
Finally, for a fixed $0<h\leq1$, $\textbf{x}\in \R^{d}$ and
$\textbf{k} = (k_1,...,k_d)    \in \Z^{d}_{+}$ set
$$g_{\textbf{k}}(\textbf{x}):=\displaystyle\prod_{j=1}^{d}g_{k_{j}}(6x_{j}/h).$$
It is easy to see that these functions give the needed partition
of unity. 
\endproof

\textbf{Proof of Lemma \ref{k5}.} Clearly the conditions of lemma yield
that whenever $|{\bf x}-{\bf a}   |>4\delta_{K}$ $$1/|\textbf{x}|_{K}\leq
\delta_{K}/|\textbf{x}|\leq\delta_{K}/(|\textbf{x}-\textbf{a}|-|\textbf{a}|)\leq\delta_{K}/(|\textbf{x}-\textbf{a}|-\delta_{K})\leq4\delta_{K}/3|\textbf{x}-\textbf{a}|.
\eqno(5)$$ It is well known that for any univariate polynomial $p$
of degree at most $n$ such that $|p|\leq 1$ in $[-a,a]$ it holds
that $|p(x)|\leq (2x/a)^{n}$ whenever $|x|>a$. Therefore using (5)
and the assumption imposed on $h_{n}$ we have
$$|h_{n}(\textbf{x})|\leq
(2|\textbf{x}-\textbf{a}|/4\delta_{K})^{2n}\leq
(2|\textbf{x}|_{K}/3)^{2n}. \eqno (6)$$ 
Now it remains to note
that by $\textbf{x}/t \in K$ it follows that $|\textbf{x}|_{K}\leq
|t|$, and thus we obtain (3) from (6). This completes the proof of
the lemma.
\endproof

 \textbf{Proof of Lemma \ref{k6}.} Recall that the support of
$g_{k}$'s consists of a pair of regular cubes with edge $h\leq 1$,
so if $A_{k}:=$supp$g_{k}$ has nonempty intersection with the unit
sphere $S^{d-1}$ then $A_{k}\subset D$, where $D$ stands for the
regular cube centered at 0 with edge 4. Let now $f_{k}$ be the
characteristic function of $A_{k}$. Since at most $2^{d}$ of
$A_{k}$'s have nonempty intersection it follows that
$$\sum f_{k}(\textbf{x})\leq 2^{d}, \textbf{x}\in \R^{d}. \eqno
(7)$$ Moreover, $m(A_{k})=2h^{d}$, where $m(.)$ stands for the
Lebesgue measure in $\R^{d}$. Using (7) we have that $$\sum
\int_{D}f_{k}dm\leq 2^{d}m(D)=8^{d}. \eqno (8)$$ Since $$\int
_{D}f_{k}dm=m(A_{k})=2h^{d}$$ whenever $A_{k}\subset D$ the
statement of the lemma easily follows from (8).
\endproof

\textbf{Proof of Lemma \ref{k3}.} Denote by $g_{k}, 1\leq k\leq N$
those functions from Lemma \ref{k4} whose support $A_{k}$ has a nonempty
intersection with $S^{d-1}$. Then by Lemma \ref{k6} $$N\leq 8^{d}/2h^{d}.
\eqno(9)$$ Moreover, by (1) $$
\displaystyle\sum_{k=1}^{N}g_{k}=1 \ \hbox{ on } \ S^{d-1}.
\eqno(10)$$

Set $$B_{k}:=A_{k}\cap S^{d-1},\quad
C_{k}:=\{\textbf{u}r(\textbf{u}): \textbf{u}\in B_{k}\}\subset
 Bd(K) , \quad 1\leq k\leq N.$$ For each $1\leq k\leq N$ choose a
point $\textbf{u}_{k}\in B_{k}$ and set
$\textbf{x}_{k}:=\textbf{u}_{k}r(\textbf{u}_{k})\in  Bd(K) .$
Furthermore let $L_{k}$ be the supporting plane to $ Bd(K) $ at the
point $\textbf{x}_{k}$ and set for $1\leq k\leq N,
L_{k}^{*}:=L_{k}\cup (-L_{k})$
$$ D_{k}:=\{\textbf{x}\in L_{k}^{*}:
\textbf{x}=t\textbf{u} \ \hbox{ for some } \ \textbf{u}\in
B_{k}, t>0\}; $$ $$ f_{k}(\textbf{x}):=g_{k}(\textbf{u}), \quad
\textbf{x}\in  Bd(K) , \textbf{x}=\textbf{u}r(\textbf{u}),\quad
\textbf{u}\in S^{d-1}$$
$$q_{k}(\textbf{x}):=g_{k}(\textbf{u}),\quad  \textbf{x}\in L_{k}^{*},   
\textbf{x}=t\textbf{u},\quad \textbf{u}\in S^{d-1},\quad
t>0.$$

Clearly, $q_{k}\in C^{\infty}(L_{k}^{*} )$ is an even positive
function which by property (2) can be extended to a regular
centrally symmetric cube $I\supset K$ so that we have on $I$
$$|\partial^{m}q_{k}/\partial x_{j}^{m}|\leq c/h^{m},\quad 1\leq j\leq d,\quad 1\leq k\leq N. \eqno (11)$$
Here and in what follows we denote by $c$ (possibly distinct)
positive constants depending only on $d, m$ and $K$. We can assume
that $I$ is sufficiently large so that
$$ I\supset G_{k}:=\{\textbf{x}\in L_{k}:
|\textbf{x}-\textbf{x}_{k}|\leq 4\delta_{K}\}, \quad 1\leq k\leq
N.$$ Then by the multivariate Jackson Theorem (see e.g. \cite{timan} ) applied to the even functions $q_{k}$ satisfying (11) for
arbitrary $m\in \N$ (to be specified below), there exist even
multivariate polynomials $p_{k}$ of total degree at most $2n$ such
that
$$||q_{k}-p_{k}||_{G_{k}^{*}}\leq c/(hn)^{m}\leq 1,\quad 1\leq k\leq N, \eqno(12)$$
where $G_{k}^{*}:=G_{k}\cup(-G_{k})$, 
$h:=n^{-\gamma}$ ($0<\gamma<1$ is specified below),
and $n$ is sufficiently large. 

We claim now that without
loss of generality it may be assumed that each $p_{k}$ is in
$H_{2n}^{d}$. Indeed, since $G_{k}^{*}\subset L_{k}^{*}$ it
follows that the homogeneous polynomial
$h_{2}:=<\textbf{x},\textbf{w}>^{2}\in H_{2}^{d}$ is identically
equal to 1 on $G_{k}^{*}$ (here $\textbf{w}$ is a properly
normalized normal vector to $L_{k}$), so multiplying the even
degree monomials of $p_{k}$ by even powers of $h_{2}$ we can
replace $p_{k}$ by a homogeneous polynomial from $H_{2n}^{d}$ so
that (12) holds. Thus we may assume that $p_{k}\in H_{2n}^{d}$ and
relations (12) hold. In particular, (12) also yields that
$$||p_{k}||_{G_{k}^{*}}\leq 2, \quad 1\leq k\leq N.\eqno(13)$$ Now consider an
arbitrary $\textbf{x}\in  Bd(K) \setminus C_{k}$. Then with some $t>1$
we have $t\textbf{x}\in L_{k}^{*}$ and $q_{k}(t\textbf{x})=0$.
Hence if $t\textbf{x}\in G_{k}^{*}$ then by (12) it follows that
$$|p_{k}(\textbf{x})|\leq |p_{k}(t\textbf{x})|\leq c/(hn)^{m}.$$
On the other hand if $t\textbf{x}\notin G_{k}^{*}$ then by (13)
and Lemma \ref{k5} we obtain
$$|p_{k}(\textbf{x})|\leq 2(2/3)^{2n}.$$ The last two estimates
yield that for every $\textbf{x}\in  Bd(K) \setminus C_{k}$ we have
$$|p_{k}(\textbf{x})|\leq c((2/3)^{2n}+(hn)^{-m}), \quad 1\leq k\leq N.\eqno(14)$$

Now let us assume that ${\bf x} \in C_k$.   

Clearly, the $C^{1+\epsilon}$ property of $ Bd(K) $ yields that
whenever $\textbf{x}\in  Bd(K) , t\textbf{x}\in L_{k}^{*}, t>1$ we
have for every $1\leq k\leq N$
$$(t-1)|\textbf{x}|=|\textbf{x}-t\textbf{x}|\leq
c\min\{|\textbf{x}-\textbf{x}_{k}|,|\textbf{x}+\textbf{x}_{k}|\}^{1+\epsilon}.
\eqno(15)$$  Obviously, for every $\textbf{u}\in B_{k}$
$$\min\{|\textbf{u}-\textbf{u}_{k}|,|\textbf{u}+\textbf{u}_{k}|\}\leq
\sqrt{d}h.$$ This and (15) yields that for $\textbf{u}\in B_{k},
\textbf{x}=\textbf{u}r(\textbf{u})\in C_{k}, t\textbf{x}\in D_{k}
(c>t>1)$ we have for $1<t<c, 0<h<c$ $$ t-1\leq
ch^{1+\epsilon},\quad D_{k}\subset G_{k}^{*}, 0<h<h_{0}.
\eqno(16)$$ 
Hence using (12), (13) and (16) we obtain for
$0<h^{1+\epsilon} \le cn^{-1}   
, 1\leq k\leq N$
$$|f_{k}(\textbf{x})-p_{k}(\textbf{x})|=|q_{k}(t\textbf{x})-p_{k}(\textbf{x})|\leq
|q_{k}-p_{k}|(t\textbf{x})+|p_{k}(t\textbf{x})-p_{k}(\textbf{x})|\leq$$
$$c/(hn)^{m}+|p_{k}(\textbf{x})|(t^{2n}-1)\leq
c(   (hn)^{-m}+nh^{1+\epsilon}),\quad \textbf{x}\in C_{k}. \eqno(17)$$

Denote for $\textbf{x}\in  Bd(K) $ $$R(\textbf{x}):= \{k:
\textbf{x}\in C_{k}\},\quad \#R(\textbf{x})\leq 2^{d}.$$

Then using the above relation together with (10),(17),(14) and (9)
we obtain for every $\textbf{x}\in  Bd(K) $
$$|1-\displaystyle\sum_{k=1}^{N}p_{k}(\textbf{x})|=|\displaystyle\sum_{k=1}^{N}(f_{k}-p_{k})(\textbf{x})|\leq
|\sum_{k\in R(\textbf{x})}...| + |\sum_{k\notin
R(\textbf{x})}...|$$ $$\leq c2^{d}(1/(hn)^{m}+nh^{1+\epsilon}) +
N||p_{k}||_{ Bd(K) \setminus C_{k}}$$ $$\leq
c( h^{-m-d}n^{-m} +h^{-d}(2/3)^{2n}+ nh^{1+\epsilon} ).\eqno(18)$$
Now it remains to choose proper values for   $m$ and
$h$. 

Choose $m \in \mathbb{N}$ to be so large that 
$$R:={m\epsilon -d \over 1+m+\epsilon+d} > \tau\epsilon \ \hbox{ and let } \ 
\gamma:= {1+m \over 1+m+\epsilon+d} . $$  
Letting $h:=n^{-\gamma}$ 
we see that $h^{-m-d}n^{-m} = nh^{1+\epsilon} = n^{-R}$.   
(Hence the $h^{1+\epsilon} \le cn^{-1} $ condition is satisfied.)  
In addition $h^{-d}(2/3)^{2n} = O(n^{-R})$, too.  
This completes the proof of Lemma \ref{k3}.
\endproof

\textbf{Proof of Theorem \ref{k1}.}  First we use the classical
Weierstrass Theorem to approximate $f\in C( Bd(K) )$ by a polynomial
$$p_{m}=\displaystyle\sum_{j=0}^{m}h_{j}^{*}, \quad h_{j}^{*}\in H_{j}^{d},\quad 0\leq j\leq m$$ of degree at most $m$ so that
$$||f-p_{m}||_{ Bd(K) }\leq \delta$$ with any given $\delta>0$.
Let $\tau \in (0,1)$ be arbitrary.  
According to Lemma \ref{k3} there exist $h_{n,j}\in
H_{2n-2[j/2]}^{d}$ such that
$||1-h_{n,j}||_{ Bd(K) }=O(n^{-\tau\epsilon}), 0\leq j\leq m$. 
Clearly,
$$h^{*}:=\displaystyle\sum_{j=0}^{m}h_{n,j}h_{j}^{*}\in
H_{2n}^{d}+H_{2n+1}^{d}$$ and $$||f-h^{*}||_{ Bd(K) }\leq \delta
+O(n^{-\tau\epsilon}).$$
\endproof

\section{Proof of Theorem \ref{k2}}











\begin{definitions}
Let $L\subset\rr$ and let $f: \ L\to\rr \cup\{-\infty \} \cup\{+\infty \} $ be a function which is defined almost everywhere (a.e.) on $L$.
We say that $f$ is {\it increasing} if $f(x)\le f(y)$ whenever $f$ is defined at $x$ and $y$ and $x\le y$. We say that $f$ is {\it increasing} almost everywhere if there exists $L^* \subset L$ such that $L\setminus L^*$ has Lebesgue measure zero, $f(x)$ is defined for all $x \in L^*$ and $f(x)\le f(y)$ whenever $x,y\in L^*, \ x\le y$.

We say that $f$ is convex if $f$ is absolutely continuous, and $f'(x)$ (which exists a.e.) is increasing a.e. on $L$.

Let $\rrr:=\rr\cup\{\infty\}$ denote the one-point compactified real line (whose topology is isomorphic to the topology of the unit circle).
\end{definitions}

Let $W: \rrr\to\rr$ be a non-negative function. Define $Q: \ \rrr\to(-\infty,+\infty]$ by 
$$W(t)=\exp(-Q(t)).$$

In the rest of the paper we will have the following assumptions on the weight $W(t)$, $t\in\rrr$:
\begin{eqnarray}\label{42}
{1\over W(t)}  \hbox{ is positive and convex on }  \rr
\end{eqnarray}
\begin{eqnarray}\label{43}
{|t|\over W(-{1\over t})} \hbox{ is positive and convex on } \rr. 
\end{eqnarray}

\begin{remark}
Equivalently, instead of \rf{43} we may assume that \rf{41} below holds and $ \lim_{t\to+\infty} $ $t(tQ'(t)-1) \le \lim_{t\to-\infty} t(tQ'(t)-1).$ We also remark that \rf{42} implies that \rf{43} is satisfied on $(-\infty,0)$ and on $(0,+\infty)$.
\end{remark}

We mention the function $W(t)=(1+|t|^m)^{-1/m}$, $1\le m$, as an example which satisfies \rf{42} and \rf{43}.

We say that a property is satisfied {\it inside} $\rr$ if it is satisfied on all compact subsets of $\rr$.

Some consequences of \rf{42} and \rf{43} are as follows.

\begin{eqnarray}\label{41}
\lim_{t\to\pm\infty} |t|W(t)=\rho\in(0,+\infty) \hbox{ exists. } 
\end{eqnarray}

Since $\exp(Q(t))$ is convex, it is Lipschitz continuous inside $\rr$. 
So $\exp(Q(t))$ is absolutely continuous inside $\rr$ which implies that both  $W(t)$ and $Q(t)$ are absolutely continuous inside $\rr$.  

$Q'(t)$ is bounded inside $\rr$ a.e. because by \rf{42} $\exp(Q(t))Q'(t)$ is increasing a.e.

We collected below some frequently used definitions and notations in the paper.
\begin{definitions}
Let $L\subset\rr$ and let $f: \ L\to\rr \cup\{-\infty \} \cup\{+\infty \} $.

$f$ is H\"older continuous with H\"older index $0<\tau \le 1$ if with some $K$ constant $|f(x)-f(y)| \le K|x-y|^\tau,$ $x,y\in L$. In this case we write $f \in H^\tau(L).$

The $L^p$ norm of $f$ is denoted by $||f||_p$. When $p=\infty$ we will also use the $||f||_{L}$ notation.

We say that an integral or limit exists if it exists as a real number.

Let $x\in\rr$. If $f$ is integrable on $ L\setminus(x-\epsilon,x+\epsilon) $ for all $0<\epsilon$ then the Cauchy principal value integral is defined as
$$ PV \int_L f(t)dt := \lim_{\epsilon\to 0^+} \int_{L\setminus(x-\epsilon,x+\epsilon)} f(t)dt ,$$
if the limit exists. 

It is known that $PV \int_L g(t)/(t-x)dt$ exists for almost every $x\in\rr$ if $g: \ L\to\rr$ is integrable. 

For $0<\iota$ and $a\in\rr$ we define
$$a_\iota^+ := \max(a,\iota) \ \hbox { and }  \ a_\iota^- := \max( -a,\iota).$$

For $a>b$ the interval $[a,b]$ is an empty set.

We say that a property is satisfied {\it inside} $L$ if it is satisfied on all compact subsets of $L$.

$o(1)$ will denote a number which is approaching to zero. For example, we may write $10^x = 100+o(1)$ as $x\to 2$. Sometimes we also specify the domain (which may change with $\epsilon$) where the equation should be considered. For example, $\sin(x) = o(1)$ for $x\in[\pi,\pi+\epsilon]$ when $\epsilon\to 0^+$.

The equilibrium measure and its support $S_w$ is defined on the next page.
Let $[a_\lambda,b_\lambda]$ denote the support $S_{W^\lambda}$ (see Lemma \rf{145}). 

For $ x\not\in(a_\lambda ,b_\lambda )$ let $V_\lambda (x):=0$, and for a.e. $x\in(a_\lambda ,b_\lambda )$ let
\begin{eqnarray}\label{90}
V_\lambda (x):={  PV \int_{a_\lambda }^{b_\lambda } { \lambda {\sqrt{(t-a_\lambda )(b_\lambda -t)} } Q'(t) \over  t-x}dt  \over \pi^2\sqrt{ (x-a_\lambda )(b_\lambda -x)} }
+{1\over \pi\sqrt{(x-a_\lambda)(b_\lambda-x)}}.
\end{eqnarray}

Let $x\in[-1,1]$. 
Depending on the value of $c\in [-1,1]$ the following integrals may or may not be principal value integrals.
$$
v_c(x):= - PV \int_{-1}^c  {\lambda \sqrt{1-t^2} e^{-Q(t)} \over \pi^2 \sqrt{1-x^2} (t-x) }dt,  
$$
\begin{eqnarray*}\label{}
h_c(x):= PV \int_c^1 {\lambda\sqrt{1-t^2} e^{-Q(t)} \over \pi^2 \sqrt{1-x^2} (t-x)}dt.
\end{eqnarray*}
(We should keep it in mind that $v_c(x)$ and $h_c(x)$ also depends on $\lambda$.)

Define
\begin{eqnarray*}\label{}
B(x)
:=v_c(x)-h_c(x)
=v_1(x)
= -PV \int_{-1}^1 { \lambda \sqrt{1-t^2} e^{-Q(t)} \over \pi^2 \sqrt{1-x^2} (t-x) }dt, \quad x\in[-1,1].    
\end{eqnarray*}

$P_n(x)$ and $p_n(x)$ denote polynomials of degree at most $n$. 
\end{definitions}

Functions with smooth integrals was introduced by Totik in \cite{totik}.
\begin{definitions}
We say that $f$ has smooth integral on $R\subset L$, if $f$ is non-negative a.e. on $R$ and
\begin{eqnarray}\label{99}
\int_I f =(1+o(1))\int_J f
\end{eqnarray}
where $I, \ J \subset R$ are any two adjacent intervals, both of which has length $0<epsilon$, and $\epsilon\to 0$. The $o(1)$ term depends on $\epsilon$ and not on $I$ and $J$.

We say that a family of functions $ {\cal F} $ has uniformly smooth integral on $R$, if any $f \in {\cal F}$ is non-negative a.e. on $R$ and \rf{99} holds, where the $o(1)$ term depends on $\epsilon$ only, and not on the choice of $f$, $I$ or $J$.
\end{definitions}

Cleary, if $f$ is continuous and it has a positive lower bound on $R$ then $f$ has smooth integral on $R$. Also, non-negative linear combinations of finitely many functions with smooth integrals on $R$ has also smooth integral on $R$.

From the Fubini Theorem it follows that if $\nu$ is a finite positive Borel measure on $T\subset \rr$ and $\{ v_t(x): \ t\in T\}$ is a family of functions with uniformly smooth integral on $R$ for which $t \to v_t(x)$ is measureable for a.e. $x \in [a,b]$, then
\begin{eqnarray*}\label{}
v(x):=\int_T v_t(x) d\nu(t)
\end{eqnarray*}
has also smooth integral on $R$.

Finally, if $f_n \to f$ uniformly a.e. on $R$, $f_n$ has smooth integral on $R$ and $f$ has positive lower bound a.e. on $R$ then $f$ has smooth integral on $R$.


\begin{remark}\label{146}
Since $\exp(-Q)$ is absolutely continuous inside $\rr$ and $(\exp(-Q))'=-\exp(-Q)Q'$ is bounded a.e. on $[-1,1]$, by the fundamental theorem of calculus we see that $\exp(-Q(t))\in H^1([-1,1])$.  
And $\sqrt{1-t}, \ \sqrt{1+t} \in H^{0.5}([-1,1])$, so $\sqrt{1-t}\sqrt{1+t}\exp(-Q(t)) \in H^{0.5}([-1,1])$ so $\sqrt{1-x^2}B(x)\in H^{0.5}([-1,1])$
by the Plemelj-Privalov Theorem (\cite{musk}, \textsection{19}). 
As a consequence, $v_c(x)$ and $h_c(x)$ exist for any $x \in [-1,1]\setminus \{c\}.$
\end{remark}

\vskip 0.5 cm

The following definitions and facts are well known in logarithmic potential theory (see \cite{st} and \cite{simeonov}).

Let $w(x)\not\equiv 0$ be a non-negative continuous function on $\rrr$ such that
\begin{eqnarray}\label{170}
\lim_{x\to \infty} |x|w(x)=\alpha \in[0,+\infty) \ \hbox{ exists }.
\end{eqnarray} 
When $\alpha=0$, then $w$ belongs to the class of so called ``admissible" weights. 

We write $w(x)=\exp(-q(x))$ and call $q(x)$ external field.
If $\mu$ is a positive Borel unit measure on $\rrr$ - in short a ``probability measure", then its weighted energy is defined by
\begin{eqnarray*}\label{}
I_w(\mu) :=\int\int \log{1 \over |x-y|w(x)w(y)} d\mu(x)d\mu(y).
\end{eqnarray*}
The integrand is bounded from below (\cite{simeonov}, pp. 3), so $I_w(\mu)$ is well defined and $-\infty <I_w(\mu)$.
Whenever it makes sense, we define the (unweighted) logarithmic energy of $\mu$ as $I_1(\mu)$ where $1$ denotes the constant 1 function.
There exists a unique probability measure $\mu_w$ - called the equilibrium measure associated with $w$ - which minimizes $I_w(\mu)$. 
Also, 
$$V_w := I_w(\mu_w) \quad \hbox{ is finite,}$$
and $\mu_w$ has finite logarithmic energy when $\alpha=0$.

If the support of $\mu$ is compact, we define its potential as
$$U^\mu(x) :=\int\log{1 \over |t-x|}d\mu(t).$$
This definition makes sense for a signed measure $\nu$, too, if $\int \Big|\log|t-x|\Big| d|\nu|(t)$ exists. 

Let
$$S_w :={\rm supp} (\mu_w) \ \hbox{ denote the support of } \ \mu_w.$$
When $\alpha=0$, then $S_w$ is a compact subset of $\rr$.
In this case with some $F_w$ constant we have
$$U^{\mu_w}+Q(x) =F_w, \quad x\in S_w.$$

\vskip 0.5 cm

Let $Bd(K)$ be the boundary of a two dimensional convex region $K \subset \rr^2$ which is centrally symmetric to the origin $(0,0)$.
For $t\in\rr$ let $(x(t),y(t))$ be any of the two points on $Bd(K) $ for which 
\begin{eqnarray}\label{177}
{y(t)\over x(t)}=t.
\end{eqnarray}
Let $x(\infty):=0$ and choose the value $y(\infty)$ such that
$(0, y(\infty) ) \in Bd(K) $. 
We define $y(\infty) / 0$ to be $\infty$, so \rf{177} also holds for $t=\infty$.

Define
\begin{eqnarray*}\label{ }
W(t):=e^{-Q(t)}:=|x(t)|, \quad t\in\rrr.
\end{eqnarray*}

\begin{lemma}
$W(t)$ satisfies properties \rf{42}, \rf{43}. And $S_W=\rrr$.
\end{lemma}
\proof
$W$ is positive on $\rr$. 

We may assume that $x(t)>0, \ t\in \rr$. Let $t_1, t_3 \in\rr$ and $t_2:=\alpha t_1+(1-\alpha)t_3$, where $0<\alpha <1$. Let $(x_2,y_2)$ be the intersection of the line segments
$\overline{ (x(t_1),y(t_1))(x(t_3),y(t_3))  }$ and $\overline{ (0,0)(x(t_2),y(t_2)) } $. Note that $1/x(t_2) \le 1/x_2$ and by elementary calculations:
$$ {1 \over x_2} 
= \alpha {1\over x(t_1)}+(1-\alpha ){1\over x(t_3)} ,
$$
so \rf{42} holds. 
The proof of \rf{43} is identical to the proof of \rf{42} once we notice that $y(-1/t)/x(-1/t) = -1/t$, and so $|t|/W(-1/t)=1/|y(-1/t)|$.

$S_W=\rrr$ follows from Corollary 3 of \cite{bdd}, since \rf{42} implies that (2.2) in \cite{bdd} is increasing on $(0,2\pi)$ with the choice $c:=0$,
and \rf{43} implies that (2.2) in \cite{bdd} is increasing on $(\pi,3\pi)$ with the choice $c:=\pi$.
(Corollary 3 can be used since \rf{41} shows that
$q(\theta ):=Q(-\cot({\theta / 2}))+\log|\sin({\theta / 2})|+\log 2$
is a continuous function on $[0,2\pi]$. And $q(\theta)$ is absolutely continuous inside $(0,2\pi)$, so it is absolutely continuous on $[0,2\pi]$.)

\endproof

\begin{lemma}\label{145}
Let $1<\lambda $. Then $S_{W^\lambda }$ is a finite interval $[a_\lambda ,b_\lambda ]$, and $ \mu_{W^\lambda } $ is absolutely continuous with respect to the Lebesgue measure and its density is $d\mu_{W^\lambda }(x) = V_\lambda (x) dx.$
\end{lemma}

\proof
Let $1<p$.
Note that $\exp(\lambda  Q(x))$ is a convex function because it is the composition of two continuous convex functions. So by \cite{b1}, Theorem 5, $S_{W^\lambda }$ is an interval $[a_\lambda ,b_\lambda ]$, which is finite since $\lim_{x\to\pm\infty} |x|W^\lambda (x)=0$. The density function $(d\mu_{W^\lambda }(x))/dx$ exists, since $(W^\lambda)'= -\exp(-\lambda Q)\lambda Q' \in L^p(\rr)$, see Theorem IV.2.2 of \cite{st}.

The integral at \rf{90} is the Hilbert transform on $\rr$ of the function defined as $\lambda \sqrt{(t-a_\lambda )(b_\lambda -t)} Q'(t) $ on $(a_\lambda ,b_\lambda )$ and $0$ elsewhere. This function is in $L^p(\rr)$, so by the M. Riesz' Theorem the integral is also in $L^p(\rr)$ hence $V_\lambda (x)$ exists for a.e. $x\in[a_\lambda ,b_\lambda ]$. Moreover, by the H\"older inequality ($1/a+1/b=1/c$ implies $||fg||_c\le ||f||_a||g||_b$) we see that $V_\lambda \in L_{1.9}(\rr)$, so $V_\lambda \in L_1(\rr)$, too.

By the proof of Lemma 16 of \cite{b0}, the function $V_\lambda $ satisfies $\int V_\lambda(x)dx =1$ and
\begin{eqnarray}\label{91}
\int_{a_\lambda }^{b_\lambda } \log|t-x|V_\lambda (t)dt = \lambda Q(x)+C, \quad x\in (a_\lambda ,b_\lambda ).
\end{eqnarray}
The left hand side is well defined since by the H\"older inequlaity 
\begin{eqnarray}\label{200}
x \ \mapsto \ \int_{a_\lambda }^{b_\lambda } \Big|\log|t-x| \Big| |V_\lambda (t)|dt
\quad \hbox{ is uniformly bounded on } \ [a_\lambda ,b_\lambda ].
\end{eqnarray}

Consider the unit signed measure $\mu$ defined by $d\mu(x):=V_\lambda (x)dx$. By \rf{91} $U^\mu(x)+\lambda Q(x)=-C$, $x\in(a_\lambda ,b_\lambda )$. 
From this and from 
$ U^{\mu_{W^\lambda} }(x) + \lambda Q(x) = F_{W^\lambda} $, $x \in[a_\lambda ,b_\lambda ]$, we get $U^\mu(x) = U^{\mu_{W^\lambda} }(x) $,  $x \in(a_\lambda ,b_\lambda )$. But \rf{200} shows that $U^{\mu^+}(x)$ and $U^{\mu^-}(x)$ are finite for all $x \in[a_\lambda ,b_\lambda ]$. So $U^{ \mu^+ }(x) = U^{ \mu_{W^\lambda} + \mu^- }(x)$, $x\in(a_\lambda ,b_\lambda )$. Here $\mu^+$ and $\mu_{W^\lambda} + \mu^-$ are positive measures which have the same mass.  $\mu_{W^\lambda} $, $\mu^-$ (and $\mu^+$) all have finite logarithmic energy (see \rf{200}), hence $ \mu_{W^\lambda}  + \mu^-$ has it, too. 
Applying Theorem II.3.2. of \cite{st} we get $U^{ \mu^+ }(z) = U^{ \mu_{W^\lambda} + \mu^- }(z) $ for all $z\in \mathbb C$. By the unicity theorem ( \cite{st}, Theorem II.2.1. ) $\mu^+  = \mu_{W^\lambda} + \mu^- $. Hence $\mu=\mu_{W^\lambda} $ and our lemma is proved.
\endproof

\begin{lemma}\label{30}
For any $[a,b]$ interval if $1<\lambda$, and $\lambda$ is sufficiently close to $1$ then $[a,b]\subset(a_\lambda,b_\lambda)$ and $V_\lambda(x)$ has positive lower bound a.e. on $[a,b]$.
\end{lemma}

\proof
First we show that $\lim_{\lambda \to 1^+} a_\lambda =-\infty$ and $\lim_{\lambda \to 1^+} b_\lambda =+\infty$. Fix $z\in\rr$ and let 
let $\lambda _n \searrow 1$ be arbitrary. We show that $z\in ( a_{\lambda _n},b_{\lambda _n})$ for large $n$. If this were not the case then for a subsequence (indexed also by $\lambda _n$) we have 
\begin{eqnarray}\label{20}
[a_{\lambda _n},b_{\lambda _n}]\subset [z,+\infty).
\end{eqnarray}
(Or, for a subsequence we have $[a_{\lambda _n},b_{\lambda _n}]\subset (-\infty,z],$ which can be handled similarly.) 
$\rrr$ is compact so by Helly's Selection Theorem (\cite{st}, Theorem 0.1.3) we can find a subsequence of the equilibrium measures $\mu_{W^{\lambda _n}}$ 
(indexed also by $\lambda_n$)
which weak-* converges to a probability measure $\mu$. This we denote by $\mu_{W^{\lambda _n}} \weak \mu$.

For fixed large $0<N$ we define the probability measure
$$\nu_N := { \mu_W \Big|_{[-N,N]} \over ||\mu_W \Big|_{[-N,N]} || } .$$
We remark that $\mu_W(\{\infty\})=0$ which implies that 
\begin{eqnarray}\label{57}
||\mu_W \Big|_{[-N,N]} || \to 1 \ \hbox{ as } N\to +\infty.
\end{eqnarray}

By (\cite{simeonov}, pp. 3) there exists $K \in \rr$ such that
\begin{eqnarray}\label{56}
K \le \log{ 1 \over |z-t| W(z)W(t) }, \quad z,t \in \rrr.
\end{eqnarray}
Now we show that 
\begin{eqnarray}\label{55}
\int\int \log{ 1 \over |z-t| W^{\lambda_1}(z)W^{\lambda_1}(t) } d\nu_N(t)d\nu_N(z) \ \hbox{ is finite. }
\end{eqnarray}
By \rf{56} the double integral at \rf{55} is bounded from below. It equals to:
\begin{eqnarray*} 
\int\int \log{ 1 \over |z-t|^{\lambda_1}  W^{\lambda_1}(z)W^{\lambda_1}(t) } d\nu_N(t)d\nu_N(z)
+ \int\int \log |z-t|^{\lambda_1-1} d\nu_N(t)d\nu_N(z).
\end{eqnarray*}
Here the first double integral is finite because $V_W$ is finite (\cite{simeonov}, Theorem 1.2). And the second integral is bounded from above since $\nu_N$ has compact support. So \rf{55} is established.

Choose $0<\tau$ such that $||\tau W(x)||_\infty\le 1$. Now,
$$
I_W(\mu) -\log(\tau^2)$$
$$
 = \lim_{M\to+\infty}\int\int \min\Big(M,\log{1\over |z-t|(\tau W(z))(\tau  W(t))}\Big)d\mu(t)d\mu(z) $$
$$
=\lim_{M\to+\infty}\lim_{n\to+\infty}\int\int \min\Big(M,\log{1\over |z-t|(\tau W(z))(\tau  W(t))}\Big)d\mu_{ W^{\lambda _n} }(t)d\mu_{W^{\lambda _n}}(z)$$
$$
\le \lim_{n\to+\infty}\int\int \log{1\over |z-t|(\tau W(z))^{\lambda _n}(\tau  W(t))^{\lambda _n}}d\mu_{W^{\lambda _n}}(t)d\mu_{W^{\lambda _n}}(z)$$
$$
\le \lim_{n\to+\infty} \int\int \log{1\over |z-t|(\tau W(z))^{\lambda _n}(\tau W(t))^{\lambda _n}}d\nu_N(t)d\nu_N(z) $$
\begin{eqnarray}\label{150}
= \int\int \log{1\over |z-t| W(z) W(t)}d\nu_N(t)d\nu_N(z) -\log(\tau^2).
\end{eqnarray}
Above in the first equality we used the monotone convergence theorem (see also \rf{56}).
In the second equality we used $\mu_{W^{\lambda _n}} \times\mu_{W^{\lambda _n}} \weak \mu \times\mu$.
In the second inequality it was used that $\mu_{W^{\lambda _n}}$ is the probability measure which minimizes the double integral of $-\log( |z-t| W^{\lambda _n}(z) W^{\lambda _n}(t))$. 
In the last equality we used the monotone convergence theorem again. (It can be used because of \rf{56}, plus the integral is finite even with the power $\lambda_1$ by \rf{55}.)

Also,
$$
\int\int \Big[ \log{1\over |z-t| W(z) W(t)} - K \Big] d\nu_N(t)d\nu_N(z) 
$$
\begin{eqnarray*} 
\le \int\int \Big[ \log{1\over |z-t| W(z) W(t)} - K \Big]
{d \mu_W(t) \over ||\mu_W \Big|_{[-N,N]} || } 
{d \mu_W(z) \over ||\mu_W \Big|_{[-N,N]} || }.
\end{eqnarray*}
Combining this with \rf{150} we have
$$I_W(\mu) 
\le K \Big[ 1 - {1 \over ||\mu_W \Big|_{[-N,N]} ||^2 }   \Big]
+ {1 \over ||\mu_W \Big|_{[-N,N]} ||^2 }  V_W.
$$
Letting $N\to+\infty$ we gain $I_W(\mu)\le V_W$.
Therefore $\mu=\mu_W$.
Thus 
$\mu_{W^{\lambda _n}}  \weak \mu_W $
which contradicts \rf{20}, since $S_W=\rrr$.

To prove the positive lower bound of $V_\lambda(x)$ a.e. on $[a,b]$, let $I:=[a-1,b+1]$. Since $W^\lambda $ is an admissible weight, we can use \cite{st}, Theorem IV.4.9., to get
\begin{eqnarray*}\label{}
\mu_{W^\lambda }\Big|_{ S_{W^{\lambda ^2}} } \ge \Big(1-{1\over \lambda ^2}\Big) \omega_{S_{W^\lambda }}\Big|_{S_{W^{\lambda ^2}}},
\end{eqnarray*}
where $\omega_{ S_{W^\lambda}}$ is the classical equilibrium measure of the set $S_{W^\lambda }$ (with no external field present). 
(We remark that $S_{W^{\lambda}} \supset S_{W^{\lambda ^2}}$.) 

It follows that if $\lambda $ is so close to $1$ that $S_{W^{\lambda ^2}} \supset I$ holds, then $[a,b]\subset(a_\lambda,b_\lambda)$ and $V_\lambda(x)$ has positive lower bound a.e. on  $[a,b]$.
\endproof

We will need Lemma 22 of \cite{b0}. We formulate it as follows:

\begin{lemma}\label{50}
Let $A<B<1$, $f\in L^1[A,1]$ and $f\in H^1[A,(B+1)/2]$. Define $v^*(x) := \int_c^1 f(t)/(t-x)dt$, where $c\in[A,B]$ and $x<c$. Then
\begin{eqnarray*}\label{}
v^*(x) = (f(c)+o(1)) \log{1 \over c-x}, \quad \hbox{ as } \quad x\to c^-.
\end{eqnarray*}
Here  $o(1)$ depends on $c-x$ only.
\end{lemma}

\begin{lemma}\label{51}
Let $-1<a<b<1$ and $0<\iota$ be fixed. Let $0<\epsilon < 1/10$ and $\delta :=\sqrt{\epsilon } - 2\epsilon $. Then for $x_1,x_2\in [a,b]\cap (c-\delta,c+\delta)^c$, $|x_1-x_2| \le \epsilon $, all the quotients
\begin{eqnarray*}\label{}
{ v_c(x_1)_\iota ^+ \over  v_c(x_2)_\iota ^+ }, \quad {v_c(x_1)_\iota ^- \over  v_c(x_2)_\iota ^-}, \quad  {h_c(x_1)_\iota ^+ \over  h_c(x_2)_\iota ^+}, \quad  {h_c(x_1)_\iota ^- \over  h_c(x_2)_\iota ^-}
\end{eqnarray*}
equal to $1+o(1)$ as $\epsilon \to 0^+$. Here the $o(1)$ term is independent of $x_1,x_2$ and $c$.
\end{lemma}
\proof
First we consider the case when $x_1, \ x_2 \le c-\delta.$ Note that for $x_1>x_2$ we have $1/(t-x_2)<1/(t-x_1)$, $t\in [c,1]$, whereas for $x_1 \le x_2$ we have
$${1 \over t-x_2} \le \Big(1+{x_2-x_1 \over c-x_2} \Big){1\over t-x_1} = (1+o(1)){1 \over t-x_1}, \quad t\in[c,1].$$
Multiplying these inequalities by $\lambda \sqrt{1-t^2} \exp(-Q(t))/\pi^2$ and integrating on $[c,1]$ we gain 
\begin{eqnarray}\label{66}
{h_c(x_2) \over h_c(x_1) } = 1+o(1),
\end{eqnarray}
where $\sqrt{1-x_2^2}/\sqrt{1-x_1^2} = 1+o(1)$ was also used.
By the same argument, if $x_1,\ x_2 \ge c+\delta$, we have $v_c(x_2)/v_c(x_1)=1+o(1)$, from which
\begin{eqnarray}\label{67}
{v_c(x_2)_\iota^+ \over v_c(x_1)_\iota^+} = 1+o(1).
\end{eqnarray}

Returning to the case of $x_1, \ x_2 \le c-\delta,$ from $v_c(x)=h_c(x)+B(x)$, from \rf{66} and from $B(x_2)=B(x_1)+o(1)$ we get
$$|v_c(x_2)-v_c(x_1)| = |o(1)| (1+|v_c(x_1)-B(x_1)|)$$
\begin{eqnarray}\label{68}
\le  |o(1)|(|v_c(x_1)|+1+||B||_{[a,b]} ).
\end{eqnarray}
Assuming $|v_c(x_1)| \le 1$, we have
\begin{eqnarray*}\label{}
|v_c(x_2)_\iota^+ - v_c(x_1)_\iota^+| \le |v_c(x_2)-v_c(x_1)| \le |o(1)|,
\end{eqnarray*}
so \rf{67} holds again.
Finally, if $|v_c(x_1)| \ge 1$, then from \rf{68}
\begin{eqnarray*}\label{}
\Big| {v_c(x_2) \over v_c(x_1)} -1 \Big| = |o(1)| \Big(1+ {1+||B||_{[a,b]}  \over |v_c(x_1)|}  \Big) = |o(1)|,
\end{eqnarray*}
from which \rf{67} again easily follows.

The proof of the rest of our lemma is similar.
\endproof

\begin{lemma}\label{79}
Let $-1<a<b<1$ and $0<\iota$ be fixed. Then the family of functions $ {\cal F}^+ :=\{ v_c(x)_\iota^+ : \ c\in [-1,1] \}$ and ${\cal F}^-:=\{ v_c(x)_\iota^- : \ c\in [-1,1] \}$ have uniformly smooth integrals on $[a,b]$.
\end{lemma}

\proof
We consider ${\cal F}^+$ only (${\cal F}^-$ can be handled similarly). Let $c\in[-1,1]$. Let $I:=[u-\epsilon ,u]$, $J:=[u,u+\epsilon ]$ be two adjacent intervals of $[a,b] $, where $0<\epsilon <1/10$. We have to show that
\begin{eqnarray*}\label{}
{ \int_I v_c(t)_\iota ^+ dt \over   \int_J v_c(t)_\iota ^+ dt } =1+o(1), \quad \hbox{ as } \quad \epsilon \to 0^+,
\end{eqnarray*}
where $o(1)$ is independent of $I,J$ and $c$.
Let $\delta :=\sqrt{\epsilon }-2\epsilon \ (>\epsilon)$.

{\it Case 1:}
Assume $I\cup J \subset (c-\delta,c+\delta)^c$. From Lemma \ref{51} we have  
$v_c(t)_\iota ^+ = (1+o(1)) v_c(t+\epsilon )_\iota ^+,$ $t\in I$. 
Thus $ \int_I v_c(t)_\iota ^+ dt = (1+o(1)) \int_J v_c(t)_\iota ^+ dt.$

{\it Case 2:}
Assume $(I\cup J) \cap (c-\delta ,c+\delta ) \not=\emptyset$. So $I\cup J\subset [ c-\sqrt{\epsilon }, c+\sqrt{\epsilon }].$
Let $\epsilon$ be so small that $c\in [(a-1)/2,(b+1)/2] $. (This can be done because of our assumption of Case 2.)

Let $f(t):=\lambda\sqrt{1-t^2} \exp(-Q(t)) / \pi^2$. 
Applying Lemma \ref{50} (with $A:=(a-1)/2$, $B:=(b+1)/2$) we have 
$\sqrt{1-x^2}h_c(x)=(f(c)+o(1)) (-\log|c-x|)$ for $x\in [c- \sqrt{\epsilon} ,c)$ as $\epsilon\to 0^+$, 
which easily leads to
\begin{eqnarray*}\label{}
h_c(x)= ({f(c) \over \sqrt{1-c^2}}+o(1)) (-\log|c-x|) \hbox{ for } x\in [c- \sqrt{\epsilon} ,c) \hbox{ as } \epsilon\to 0^+.
\end{eqnarray*}
From here using $h_c(x)=v_c(x)-B(x)$ we get
\begin{eqnarray}\label{63}
v_c(x)= ( {f(c) \over \sqrt{1-c^2}} +o(1)) (-\log|c-x|) \hbox{ for } x\in [c-\sqrt{\epsilon} ,c) \hbox{ as } \epsilon\to 0^+.
\end{eqnarray}
Clearly, \rf{63} also holds for $x\in (c,c+\sqrt{\epsilon} ]$ (which can be seen by stating Lemma \ref{50} for $-1<A<B$ instead of $A<B<1$). 

$f(x)$ has a positive lower bound on $[(a-1)/2,(b+1)/2]$. So we can choose $\epsilon $ so small that the right hand side of \rf{63} is at least $\iota $ for all possible values of $c$ and $x$. Hence $v_c(x)=v_c(x)_\iota ^+$ and
\begin{eqnarray*}\label{}
{ \int_I v_c(t)_\iota ^+ dt \over   \int_J v_c(t)_\iota ^+ dt }  
=  { ( {f(c) \over \sqrt{1-c^2}} +o(1)) \int_I \log{1\over  |c-t|} dt  \over  ( {f(c) \over \sqrt{1-c^2}} +o(1)) \int_J \log{1\over  |c-t|} dt     }
=(1+o(1))^2
=1+o(1),
\end{eqnarray*}
where we used that $\log(1/|x|)$ has smooth integral on $[-2,2] $ (\cite{b0}, Proposition 20).
\endproof

\begin{lemma}\label{78}
Let $F(x)=G(x)-H(x)$, where $F(x), \ G(x), \ H(x)$ are a.e. non-negative functions defined on an interval, $G(x)$ and $H(x)$ have smooth integrals and $H(x)\le (1-\eta)G(x)$ a.e. with some $\eta\in(0,1)$. Assume also that $\int_I F=0$ implies $\int_I G=\int_I H =0$,  when the interval $I$ is small enough. Then $F(x)$ has smooth integral.
\end{lemma}

\proof
Let $I$ and $J$ be two adjacent intervals of equal lengths $\epsilon$, where $\epsilon$ is ``small enough". Let $a:=\int_I G$, $A:=\int_J G$, $b:=\int_I H$, $B:=\int_J H$. By assumption 
\begin{eqnarray}\label{71}
A=(1+o(1))a \quad \hbox{ and } \quad B=(1+o(1))b, \quad \hbox{ as } \epsilon\to 0^+
\end{eqnarray}
and we have to show that $A-B=(1+o(1))(a-b).$ 

We may assume that $a-b\not= 0,$ otherwise $a=b=0$ from the assumption of the lemma and so $A=B=0.$

Integrating $H\le (1-\eta)G $ on $I$ we get $b\le (1-\eta)a,$ from which $(a+b)/$ $(a-b)\le (1+(1-\eta))/(1-(1-\eta)).$ Thus, from \rf{71}
\begin{eqnarray*}\label{}
|(A-a)-(B-b)| \le |o(1)|(a+b)\le |o(1)|(a-b). 
\end{eqnarray*}
\endproof

Following the proof of Lemma 24 of \cite{b0} we will prove the following lemma. 
But we remark that the absolutely continuous hypothesis of Lemma 24 is unnecessary at \cite{b0}.

\begin{lemma}\label{140}
Let $N(x)$ be a bounded, increasing, right-continuous function on $[-1,1]$ and let $f(x) \in L^1([-1,1])$ be non-negative. Then
\begin{eqnarray}\label{110}
PV \int_{-1}^1 {f(t)N(t)\over  t-x}dt = - N(1)f_{1}(x) + \int_{(-1,1]} f_t(x)dN(t), \quad a.e. \ x\in[-1,1],
\end{eqnarray}
where the integral on the right hand side is a Lebesgue-Stieltjes integral and
\begin{eqnarray*}\label{}
f_c(x):= -PV \int_{-1}^c {f(t)\over  t-x}dt, \quad a.e. \ x\in[-1,1].
\end{eqnarray*}
\end{lemma}

\proof
Let us denote the left hand side of \rf{110} by $F(x)$.
Since $f(x)$ and $f(x)N(x)$ are in $L^1[-1,1]$ and $N(x)$ is increasing, there is a set of full measure in $(-1,1)$ where $f_1(x)$, $F(x)$ and $N'(x)$ all exist. Let $x$ be chosen from this set. It follows that $f_c(x)$ exist for all $c\in[-1,1]\setminus\{x\}$. Also,
\begin{eqnarray}\label{178}
F(x) = \lim_{\epsilon  \to 0^+} \Big( \int_ {-1}^{x-\epsilon } {f(t)N(t)\over t-x}dt + \int_{x+\epsilon  }^1 {f(t)N(t)\over t-x}dt \Big).
\end{eqnarray}
$t \to f_t(x)$ is a continuous increasing function on $[-1,x)$ and it is a continuous decreasing function on $(x,1]$ so at \rf{178} we can use integration by parts to get
$$
 \int_{-1}^{x-\epsilon} + \int_{x+\epsilon  }^1 
= -f_{x-\epsilon }(x)N(x-\epsilon ) + f_{-1}(x)N(-1) + \int_{(-1,x-\epsilon ]}  f_t(x)dN(t)
$$
\begin{eqnarray*}\label{}
 +f_{x+\epsilon }(x)N(x+\epsilon ) - f_{1}(x)N(1) + \int_{(x+\epsilon ,1]} f_t(x)dN(t)
\end{eqnarray*}
But above $f_{-1}(x)=0$ and 
$$
f_{x+\epsilon }(x)N(x+\epsilon ) - f_{x-\epsilon }(x)N(x-\epsilon ) 
$$
\begin{eqnarray}\label{120}
= [f_{x+\epsilon }(x) - f_{x-\epsilon }(x)]N(x+\epsilon ) 
+f_{x-\epsilon }(x) [N(x+\epsilon ) -N(x-\epsilon ) ].
\end{eqnarray}
Note that
\begin{eqnarray*}\label{}
f_{x+\epsilon }(x) - f_{x-\epsilon }(x) = -PV \int_{x-\epsilon }^{x+\epsilon } {f(t) \over t-x}dt\to 0 \quad \hbox{ as } \epsilon\to 0^+,
\end{eqnarray*}
since $f_1(x)$ exists. Also, $0 \le f_{x-\epsilon }(x) \le c_1\log(1/\epsilon )$, $0<\epsilon<1$, which implies that the second term at \rf{120} tends to zero (since $N$ is differentiable at $x$).

Putting these together, we get that on one hand,
\begin{eqnarray}\label{121}
\lim_{\epsilon \to 0^+} \Big( \int_{(-1,x-\epsilon ]}   f_t(x)dN(t) + \int_{(x+\epsilon,1]} f_t(x)dN(t) \Big)
\end{eqnarray}
exists and equals to $F(x)+f_{1}(x)N(1)$, and on the other hand, \rf{121} equals to 
\begin{eqnarray}\label{122}
\int_{   (-1,1]\setminus\{x\}      } f_t(x)dN(t) = \int_{(-1,1]} f_t(x)dN(t)
\end{eqnarray}
by the monotone convergence theorem 
(which can be used since $c \to f_c(x)$ is bounded from below on $[-1,1]$ since $f_1(x)$ is finite).
The the continuity of $N$ at $x$ allowed us to integrate on the whole $[-1,1]$ at \rf{122}.
\endproof

\begin{lemma}\label{127}
Let $[a,b]$ be arbitrary and let $1<\lambda $ be chosen to satisfy the conclusion of Lemma \ref{30}. 
Then $V_\lambda(x)$ has smooth integral on $[a,b]$.
\end{lemma}

\proof
To keep the notations simple we will assume that $-1<a<b<1$, and $a_\lambda =-1$, $b_\lambda =1$, that is, the support of $\mu_{W^\lambda }$ is $[-1,1]$.
This can be done without loss of generality.
Define 
\begin{eqnarray*}\label{}
v(t) := {\lambda \sqrt{1-t^2}e^{-Q(t)}\over \pi^2\sqrt{1-x^2}} \ \hbox{ and } \ M(t):=\lim_{s\to t^+} e^{Q(s)}Q'(s),
\end{eqnarray*}
where $v(t)$ also depends on the choice of $x$. 
Note that $M(t)$, $t\in[-1,1]$, is a bounded, increasing, right-continuous function which agrees with $\exp(Q(t))Q'(t)$ almost everywhere.

Applying Lemma \ref{140} for $f(t):=v(t)$ and $N(t):=M(t)$,
let us fix an $x \in [a,b]$ value for which both \rf{110} and  $d\mu_{W^\lambda }(x) = V_\lambda (x) dx$ are satisfied. (These are satisfied almost everywhere.) 
From \rf{90} and Lemma \ref{140} we have
$$
V_\lambda (x)
= {1 \over \pi\sqrt{1-x^2}} + PV \int_{-1}^1 {\lambda \sqrt{1-t^2}Q'(t) \over \pi^2\sqrt{1-x^2}(t-x)}dt 
$$
\begin{eqnarray*}\label{}
= {1 \over \pi\sqrt{1-x^2}} + PV \int_{-1}^1 {v(t)M(t)\over  t-x}dt
= L(x) + \int_{(-1,1 ]}   v_t(x)dM(t),
\end{eqnarray*}
where $L(x):= 1 / (\pi\sqrt{1-x^2}) - M(1)B(x) $. 

Let $0<\iota$.
Since $L(x)$ is a continuous function on $[a,b]$ (see Remark \ref{146}), $L(x)_\iota^+$ and $L(x)_\iota^-$ have smooth integrals on $[a,b]$.
Also, by Lemma \ref{79} ${\cal F}^+$ and  ${\cal F}^-$ have uniformly smooth integrals on $[a,b]$, so both
$$
V_\lambda (x)_{(\iota)}^{(+)} := L(x)_\iota^+ + \int_{(-1,1 ]}   v_t(x)_\iota^+ dM(t) \quad \hbox{ and } 
$$
\begin{eqnarray*}\label{}
V_\lambda (x)_{(\iota)}^{(-)} := L(x)_\iota^- + \int_{(-1,1 ]}   v_t(x)_\iota^- dM(t) 
\end{eqnarray*}
have smooth integrals on $[a,b]$. (These new functions are not to be mixed with $V_\lambda (x)_{\iota}^{-}$ and $V_\lambda (x)_{\iota}^{-}$.)

Set
\begin{eqnarray*}\label{}
V_\lambda (x)_{(\iota)} := V_\lambda (x)_{(\iota)}^{(+)} - V_\lambda (x)_{(\iota)}^{(-)}.
\end{eqnarray*}
Then, using $| z_\iota^+ - z_\iota^- - z| \le \iota, \ z\in\rr$, we get
$$
| V_\lambda(x)_{(\iota)} - V_\lambda (x)| 
\le |L(x)_\iota^+ - L(x)_\iota^- - L(x)| +  
 \int_{(-1,1]}   |v_t(x)_\iota^+ - v_t(x)_\iota^- - v_t(x)| dM(t)
$$
\begin{eqnarray}\label{80}
\le \iota + \int_{(-1,1]}   \iota dM(t)
= \iota (1+M(1)-M(-1)).
\end{eqnarray}
So 
\begin{eqnarray}\label{105}
V_\lambda(x)_{(\iota)} \to V_\lambda (x) \ \hbox{ uniformly a.e. on } [a,b] \hbox{ as } \iota \to 0^+.
\end{eqnarray}

And since
\begin{eqnarray}\label{81}
V_\lambda (x) \ \hbox{ has positive lower bound a.e. on } \ [a,b],
\end{eqnarray}
$V_\lambda (x)_{(\iota)}$ has also positive lower bound a.e. on $[a,b]$, assuming $\iota$ is small enough.
In addition, $v_t(x)\ge 0$ when  $t\in[0,x]$, whereas $v_t(x)\ge B(x)\ge -||B||_{[a,b]} $ when $t\in(x,1]$,
so $V_\lambda (x)_{(\iota)}^{(-)}$ is bounded a.e. on $[a,b]$. 
It follows that $V_\lambda (x)_{(\iota)}^{(-)} \le (1-\eta)V_\lambda(x)_{(\iota)}^{(+)}$ a.e. $x\in [a,b]$ for some $\eta\in(0,1)$.

Applying Lemma \ref{78} we conclude that $V_\lambda (x)_{(\iota)}$ has smooth integral on $[a,b]$ (if $\iota$ is small enough). 
Therefore $V_\lambda (x)$ has smooth integral by \rf{105} and \rf{81}.
\endproof

Approximation by weighted polynomials with varying weights was introduced by Saff (\cite{saff}). In our proof we shall utilize the strong connection between weighted polynomials and homogeneous polynomials on the plane.

It was proved by Kuijlaars (\cite{kuijlaars}, see also \cite{st}, Theorem VI.1.1) that when $\alpha=0$ at \rf{170} then there exists a closed set $Z(w)\subset\rr$ with the property that a continuous function $f(x), \ x\in\rr$, is the uniform limit of weighted polynomials $w^nP_n$ $(n=0,1,2,...)$ on $\rr$ if and only if $f(x)$ vanishes on $Z(w)$.
We formulate the following version of this theorem.

\begin{lemma}\label{39}
Assume that $0<\alpha$ at \rf{170}. 
Then there exists a closed set $ Z_\rrr(w) $ such that a continuous function $f(x), \ x\in\rrr$, is the uniform limit of weighted polynomials $ w^n p_n$ $(n=0,2,4,...)$ on $\rrr$ if and only if $f(x)$ vanishes on $Z_\rrr(w)$.
\end{lemma}

\proof
Let $X:=\rrr$. 
Note that $w^n p_n$ is continuous on $\rrr$ when $n$ is even. (Naturally the value $(w^n p_n)(\infty)$ is defined to be $\lim_{x\to \pm \infty} (w^n p_n)(x)$.)

Let ${\cal A}$ be the collection of continuous functions $f$ on $X$ such that $w^n p_n \to f$ $(n=0,2,4,...)$ uniformly on $X$ for some $p_n$. 
Define the set $ Z_\rrr(w)  :=\{x\in X: \ f(x)=0 \hbox{ for all } f\in {\cal A} \}$, which is certainly closed.

It is easy to see (similarly as in \cite{st}, Theorem VI.1.1) that  ${\cal A}$ is an algebra which is closed under uniform limits. Also, it separates points in the sense that if $x_1, x_2 \in X\setminus  Z_\rrr(w) $ are two distinct points, then there exists $f \in {\cal A} $ such that $f(x_1)\not= f(x_2)$. Indeed, let us assume that, say, $x_2$ is finite and let $g \in {\cal A}$ such that $g(x_1)\not=0$. Let $w^n p_n \to g$ $(n=0,2,4,...)$ uniformly on $X$. 
Then $w^{n+2}(x)[(x-x_2)^2 p_n(x)] \to w^2(x)(x-x_2)^2g(x)=:f(x)$ $(n=0,2,4,...)$ uniformly on $X$ 
because $||w^2(x)(x-x_2)^2||_\rrr < +\infty$. Thus $f(x) \in{\cal A} $. 
And $f(x_1)\not=0=f(x_2)$ (which holds even if $x_1$ was infinity.)

Since ${\cal A}$ satisfies the properties above, by the Stone-Weierstrass Theorem
\begin{eqnarray*}\label{}
{\cal A} = \{f: \ f \hbox{ is continuous on } X \hbox{ and } f\equiv 0 \hbox{ on }  Z_\rrr(w)  \}.
\end{eqnarray*}
\endproof


We now restate Theorem \ref{k10} and prove it.

\begin{theorem}\label{201}
For a weight satisfying \rf{42} and \rf{43} we have $ Z_\rrr(W)  =\emptyset.$ 
That is, any continuous function $g:\rrr\to\rr$ can be uniformly approximated by weighted polynomials $W^n p_n$ $(n=0,2,4,...)$ on $\rrr$.
\end{theorem}
\proof
Let $x_0\in\rrr$. We show that $x_0\not\in Z_\rrr(W)$. 

First let us assume that $x_0$ is finite.
Choose $J:=[a,b]$ such that $a<x_0<b$ holds. 
Let $f(x)$ be a continuous function which is zero outside $J$ and $f(x_0)\not= 0$. 
Let $1<\lambda =u/v$ ($u,v\in\mathbb N^+$) be a rational number for which the conclusion of Lemma \ref{30} holds. 
Now we use a powerful theorem of Totik.
Since $V_\lambda $ has a positive lower bound a.e. on $J$ and it has smooth integral on $J$ (see Lemma \ref{127}), by \cite{totik}, Theorem 1.2, $(a,b)\cap {Z}(W^\lambda )=\emptyset$. 
So we can find $P_n$ $(n=0,1,2,...)$ such that $(W^\lambda )^nP_n \to f$ uniformly on $\rrr$.

So for $n:=Nv,$ we have
\begin{eqnarray}\label{31}
W^{Nu} p_{Nu} \to f, \quad N=0,1,2,..., \ \hbox{ uniformly on } \ \rrr, 
\end{eqnarray}
where $p_{Nu}:=P_{Nv}$ and $\deg (p_{Nu})\le Nv \le Nu$.
For all fixed $s\in\{0,...,u-1 \}$ if we approximate $f/W^s$ instead of $f$ at \rf{31}, it easily follows that 
there exist $p_k$ $(k=0,1,2,...)$ such that 
\begin{eqnarray}\label{38}
W^{k} p_{k} \to f, \quad k=0,1,2,..., \ \hbox{ uniformly on } \ \rrr. 
\end{eqnarray}
Using only $k=0,2,4,...$, we get $x_0\not\in Z_\rrr(W)$ by Lemma \ref{39}.

Now let $x_0=\infty$. Define
$$W_0(x):={1 \over |x|} W(-{1 \over x}).$$
Note that $1/W_0(x) \ (= |x|/W(-1/x) )$ and $|x|/W_0(-1/x) \ (= 1/W(x))$ are positive and convex functions because $W$ satisfies \rf{43} and \rf{42}. 

Let $g$ be a continuous function on $\rrr$. Define $-1/\infty$ to be $0$ and $-1/0$ to be $\infty$. (So $g(x)$ is continuous on $\rrr$ if and only if  $g(-1/x)$ is continuous on $\rrr$.)
Observe that for some $p_n$ we have 

\noindent
$W^n(x)p_n(x) \to g(x)$ $(n=0,2,4,...)$ uniformly on $\rrr$, iff 

\noindent
$W^n(-1/x)p_n(-1/x) \to g(-1/x)$ $(n=0,2,4,...)$ uniformly on $\rrr$, iff 

\noindent
${W_0}^n(x) q_n(x), \to g(-1/x)$ $(n=0,2,4,...)$ uniformly on $\rrr$, 

\noindent
where $q_n(x):=x^n p_n(-1/x)$ are polynomials, $\deg q_n \le n$. 

Now let $f(x)$ be a continuous function on $\rrr$ which is zero in a neighborhood of $0$ but $f(\infty)\not=0$. 
By what we have already proved, $q_n$ polynomials exist such that ${W_0}^n(x)q_n(x)$ $(n=0,2,4,...)$ tends to $f(-1/x)$ uniformly.
Therefore we can approximate $f(x)$ uniformly by $W^n(x)p_n(x)$ $(n=0,2,4,...)$, where $p_n(x):=x^n q_n(-1/x)$.

\endproof

\begin{lemma}\label{157}
Let $f(x,y), \ (x,y) \in Bd(K) $, be a continuous function such that $f(x,y)=f(-x,-y)$ for all $(x,y)\in Bd(K) $. Then homogeneous polynomials
$$ h_n(x,y) := \sum_{k=0}^n a_k^{(n)} x^{n-k}y^k, \quad n=0,2,4,...$$
exist such that $h_n(x,y) \to f(x,y)$ $(n=0,2,4,...)$ uniformly on $ Bd(K) $.
\end{lemma}
\proof
Recall the definition: $y(t)/x(t)=t, \ t\in\rrr,$ where $(x(t),y(t))\in Bd(K) $ and $W(t):= |x(t)|$.
Define
$$f(t):=f(x(t),y(t))=f(-x(t),-y(t)), \quad t\in\rrr.$$
Note that if $n$ is an even number (and $a_k^{(n)}$ are unknowns) then
$$\sum_{k=0}^n a_k^{(n)}  x^{n-k}(t) y^k(t)
= x^n(t) \sum_{k=0}^n a_k^{(n)}  \Big({y(t) \over x(t)} \Big)^k$$
\begin{eqnarray}\label{202}
= |x(t)|^n \sum_{k=0}^n a_k^{(n)}  t^k
=W^n(t) p_n(t),
\end{eqnarray}
where $p_n(t):=\sum_{k=0}^n a_k^{(n)}  t^k$, $\deg p_n \le n$.
(When $t=\infty$, the left hand side of \rf{202} again equals to $(W^n p_n)(\infty) := \lim_{t\to \pm \infty} W^n(t)p_n(t)$.)

But by Theorem \ref{201} there exist $W^n(t)p_n(t)$ $(n=0,2,4,...)$ which tends to $f(t)$ uniformly on $\rrr$. This completes the proof, since for any $(x,y) \in  Bd(K) $ there exists $t\in\rrr$ such that either $(x(t),y(t))=(x,y)$ or $(-x(t),-y(t))=(x,y)$.
\endproof

\noindent
{\bf Proof of Theorem \ref{k2}.}

\noindent
Define $f(x,y):=1, \ (x,y) \in Bd(K)$. By Lemma \ref{157} there exist $h_{2n} \in H_{2n}^2$, $n \in \mathbb{N}$, such that $||1-h_{2n}||_{ Bd(K) } \to 0$. From here Theorem \ref{k2} follows the same way Theorem \ref{k1} follows from Lemma \ref{k3}.
\endproof



\

\noindent
D. Benko \\
Department of Mathematics \\
Western Kentucky University \\
Bowling Green, KY 42101 \\
USA \\
E-mail: {\it dbenko2005@yahoo.com} \\

\noindent
A. Kro\'o \\
Alfr\'ed R\'enyi Institute of Mathematics \\
Hungarian Academy of Sciences \\
H-1053 Budapest, Re\'altanoda u. 13-15 \\
Hungary \\
E-mail: {\it kroo@renyi.hu} \\

\end{document}